\documentclass[12pt]{amsart}
\usepackage{amsthm}
\usepackage{amsmath, mathrsfs}
\usepackage{amssymb}
\usepackage{hyperref}
\usepackage{enumerate}
\usepackage{latexsym,array}
\usepackage{amsfonts}
\usepackage{shadow}

\newcommand{\bigH}{\mathcal{H}}

\newcommand{\bigS}{\mathcal{S}}

\newcommand{\bigW}{\mathcal{W}}

\newcommand{\bigK}{\mathcal{K}}

\newtheorem{Pa}{Paper}[section]
\newtheorem{Tm}[Pa]{{\bf Theorem}}

\newtheorem{Dn}[Pa]{{\bf Definition}}

\newtheorem{Rk}[Pa]{{\bf Remark}}
\newtheorem{Pn}[Pa]{{\bf Proposition}}

\date{}
\keywords{convolution algebra, non-commutative white noise space,
non-commutative stochastic distributions} \subjclass{Primary:
16S99, 60H40, 93B07. Secondary: 93A25}
\thanks{D. Alpay thanks the
Earl Katz family for endowing the chair which supported his
research, and the Binational Science Foundation Grant number
2010117.
}
\author{Daniel Alpay}
\author{Guy Salomon}
\address{Department of Mathematics
\newline
Ben Gurion University of the Negev \newline P.O.B. 653,
\newline
Be'er Sheva 84105, \newline ISRAEL} \email{(DA)
dany@math.bgu.ac.il} \email{(GS) guysal@math.bgu.ac.il}
\title[Non-commutative stochastic distributions]
{Non-commutative stochastic distributions and applications to
linear systems theory}
\begin{document}


\maketitle


\begin{abstract}
In this paper, we introduce a non-commutative space of stochastic
distributions, which contains the non-commutative white noise
space, and forms, together with a natural multiplication, a
topological algebra. Special inequalities which hold in this
space allow to characterize its invertible elements and to
develop an appropriate framework of non-commutative stochastic
linear systems.
\end{abstract}

\section{Introduction}
\setcounter{equation}{0}
In this paper we introduce and study a non-commutative version of
a space of stochastic distributions, and give applications to
mathematical system theory.  To set the problem into perspective,
recall that, in white noise analysis, various spaces of
stochastic distributions have been introduced by Hida, Kondratiev,
and others; see \cite{MR1408433} and the references therein. Among
those introduced by Kondratiev, one (denoted by $\mathcal
S_{-1}$) plays an important role. It is the dual of a Fr\'echet
nuclear space, and in particular the increasing union of a
countable family of Hilbert spaces with decreasing norms.
$\mathcal S_{-1}$ is an algebra when endowed with the Wick
product, and the Wick product satisfies in $\mathcal S_{-1}$ an
inequality, called V\r{a}ge inequality. The space $\mathcal
S_{-1}$ was recently used to develop a new approach to the theory
of linear stochastic systems, when not only the input is random
but also the characteristics of the system. See
\cite{aa_goh,alp,MR2610579}. We recently defined a large class of
topological algebras, which also satisfy a V\r{a}ge type
inequality, and which are furthermore closed under tensor
products. See \cite{vage1,2012arXiv1204.5277A}. For the non-commutative
version of the white noise and of the white noise
space we refer to \cite{MR1217253}. The non-commutative
counterparts of spaces of stochastic distributions, especially
ones which satisfy V\r{a}ge type inequalities, do not seem to
have been studied. We begin such a study here, and give
applications to non-commutative linear systems parallel to the
one done in \cite{aa_goh,alp,MR2610579} for the Kondratiev space
and in
\cite{vage1} for V\r{a}ge spaces.\\

We divide this introduction into three parts. The first two parts
are preliminaries about the commutative case, namely on the white
noise space and on the Kondratiev space $\mathcal S_{-1}$ of
stochastic distributions. In the third part we discuss our
approach to define a non-commutative space of stochastic
distributions and
give an outline of the paper.\\

\subsection{The (commutative) white noise space}
To set the framework of the commutative case we recall the
following definitions. Let $\bigH$ be a separable complex Hilbert space. We
consider its $n$-fold Hilbert spaces tensor power $\bigH^{\otimes
n}$. The symmetric product $\circ$ is defined by
\[
u_1 \circ \cdots \circ u_n=\frac{1}{n!} \sum_{\sigma \in S_n}
u_{\sigma(1)}\otimes \cdots \otimes u_{\sigma(n)},
\]
and the closed subspace of $\bigH^{\otimes n}$ generated by all
vectors of this form is called the $n$-th symmetric power of
$\bigH$, and denoted by $\bigH^{\circ n}$. See \cite{Neveu68}. We
make the convention $\bigH^{\otimes 0}=\mathbb C$, and the
element $1 \in \mathbb C$ is called the vacuum vector and denoted
by $\bf 1$. Two inner products are defined on $\bigH^{\circ n}$.
The first is called the symmetric inner product, and defined by
\[
\langle u_1 \circ \cdots \circ u_n , v_1 \circ \cdots \circ v_n
\rangle_\circ = per(\langle u_i,v_j\rangle),
\]
where $per(A)$ is called the permanent of $A$ and has the same
definition as a determinant, with the exception that the factor
$sgn(\sigma)$ is omitted. The second is called the tensor inner
product. It is induced by the tensor inner product on
$\bigH^{\otimes n}$
\[
\langle u_1 \otimes \cdots \otimes u_n , v_1 \otimes \cdots
\otimes v_n \rangle_{\otimes} =\prod_{i=1}^n\langle u_i,v_i\rangle.
\]
Therefore, the tensor inner product on $\bigH^{\circ n}$ is simply
\[
\langle u_1 \circ \cdots \circ u_n , v_1 \circ \cdots \circ v_n
\rangle_{\otimes}=\frac{1}{n!^2}\sum_{\sigma, \tau \in S_n}
\langle u_{\sigma(1)},v_{\tau(1)} \rangle \cdots \langle
u_{\sigma(n)},v_{\tau(n)} \rangle.
\]
It is clear that $\|\cdot\|_\otimes=\frac1{n!}\|\cdot\|_\circ$.
Assuming $(e_i)_{i \in I}$ is an orthonormal basis of $\bigH$ where $I\subseteq \mathbb N$,
for $\alpha:I \to \mathbb N_0$ (for simplicity, we denote $\alpha_i$ instead of $\alpha(i)$) with a support
$\{i_1,\ldots,i_m\}$ ($i_1<\cdots<i_m$) such that
$|\alpha|=\sum_{j=1}^m \alpha_{i_j}=n$, we denote
\[
e_\alpha=e_{i_1}^{\circ \alpha_{i_1}}\circ \cdots \circ
e_{i_m}^{\circ \alpha_{i_m}} \in \bigH^{\circ n}.
\]
$(e_\alpha)$ is clearly an orthogonal basis of $\bigH^{\circ n}$.
The squared symmetric norm of $e_{\alpha}$ is $\alpha!=\alpha_{i_1}!\alpha_{i_2}!\cdots\alpha_{i_m}!$, and
the squared tensor norm is $\frac{\alpha!}{n!}$.\\

The symmetric Fock space over $\bigH$ is the Hilbert space
\[
\Gamma^\circ(\bigH)=\oplus_{n=0}^\infty\bigH^{\circ n},
\]
with the corresponding symmetric inner product.\\

For the definition
of the white noise space, one usually takes $\bigH=\mathbf L^2(\mathbb
R)$. Let $(e_n)_{n\in\mathbb N}$ be an orthonormal basis of
$\mathbf L^2(\mathbb R)$ (for example, the Hermite functions). We
define the (commutative) white noise space  $\bigW$ as the
symmetric Fock space of $\bigH=\mathbf L^2(\mathbb R)$. Thus,
denoting by $\ell$ the free commutative monoid generated by
$\mathbb N_0$, that is,
\[
\ell= \mathbb N_0^{(\mathbb N)}=\left\{\alpha \in \mathbb
N_0^{\mathbb N} : \text{supp}(\alpha) \text{ is finite} \right\},
\]
and setting $\nu(\alpha)=\alpha!$ we conclude that
\[
\bigW=\Gamma^\circ(\bigH)=\left\{\sum_{\alpha \in \ell}f_\alpha
e_\alpha : \sum_{\alpha \in \ell}|f_\alpha|^2 \alpha!<\infty
\right\} ={\mathbf L}^2(\ell, \nu).
\]

For more information on symmetric and non-symmetric Fock spaces
we refer to \cite{MR1222649,Neveu68}.\\

In this paper, we do not use any realization of the white noise space. Nevertheless, it is worth to mention that the classical realization is as the ${\mathbf L}^2$-space of Gaussian white noise.
More precisely, given a nuclear countably Hilbert space $E$ which is densely and continuously imbedded in
$\mathbf L^2(\mathbb R)$, the Bochner-Minlos theorem insures the existence of a probability measure $P$ on the Borel
$\sigma$-algebra of $E^\prime$
such that $e^{-\frac12 \|\varphi\|^2_{{\mathbf L}_2}}
=\int_{E^\prime}e^{i\langle f, \varphi \rangle}dP(f)$.
The space ${\mathbf L}_2(E^\prime,\mathcal B, P)$ is called the Gaussian white noise space, and it is isomorphic to the
symmetric Fock space $\Gamma^\circ(\bigH)$, via the Wiener-It\^o-Segal isomorphism. For more information,
see for instance \cite[pp 162-163]{Huang}.

\subsection{The Wick product and the (commutative) Kondratiev space of
stochastic distributions}

The standard multiplication of two elements in the white noise
space is called the Wick product.
\begin{Dn}
The Wick product is defined by $(f,g) \mapsto f \circ g$ whenever
it make sense. In terms of the basis, we obtain that
\[
f\circ g = \left(\sum_{\alpha \in \ell}f_\alpha e_\alpha
\right)\circ \left(\sum_{\alpha \in \ell}g_\alpha e_\alpha
\right)=\sum_{\alpha \in \ell} \left( \sum_{\beta \leq
\alpha}f_\beta g_{\alpha-\beta} \right) e_\alpha.
\]
\end{Dn}
As it is obvious from its definition, the Wick product is
actually a convolution of functions over the monoid $\ell$. It is
well known that $\mathcal W$ is not closed under it; see Remark
\ref{notclosed}. On the other hand, the dual of the Kondratiev
space $\bigS_1$ of stochastic test functions, namely the
Kondratiev space $\bigS_{-1}$ of stochastic distributions, is
closed under the Wick product. The space $\bigS_{1}$ is defined
as follows:
\[
\bigS_{1}=\left\{ \sum_{\alpha \in \ell} f_\alpha e_\alpha:
\sum_{\alpha \in \ell}
 |f_{\alpha}|^2(2\mathbb N)^{\alpha p}(\alpha!)^2< \infty \text { for all } p
 \in \mathbb N \right\},
\]
where $(2\mathbb N)^\alpha=2^{\alpha_1} \cdot 4^{\alpha_2}\cdot
6^{\alpha_3} \cdots$. It is a countably normed Hilbert space (in
the language of Gelfand) which is a subspace of the white noise
space $\bigW$. Its dual with respect to the center space $\bigW$, namely, the Kondratiev space of stochastic
distributions $\bigS_{-1}$, can be viewed as
\[
\begin{split}
\bigS_{-1}&=\left\{ \sum_{\alpha \in \ell} f_\alpha e_\alpha:
\sum_{\alpha \in \ell}
 |f_{\alpha}|^2(2\mathbb N)^{-\alpha p}< \infty \text { for some } p
 \in \mathbb N \right\}\\
&= \bigcup_{p} {\mathbf L}^2(\ell,\mu_{-p}),
\end{split}
\]
where $\mu_{-p}$ is the point measure defined by
\[
\mu_{-p}(\alpha)=(2\mathbb N)^{-\alpha p}.
\]
Together with the white noise space these two spaces form the
Gelfand triple $(\bigS_1,\mathcal W, \bigS_{-1})$. These two
spaces $\bigS_1$ and $\bigS_{-1}$ are both nuclear (the latter
when endowed with the strong topology), a property which allows to
consider $\rm{Hom}(\bigS_1,\bigS_{-1})$ as an appropriate
framework for the theory of stochastic linear systems thanks to
Schwartz' kernel theorem; see \cite{yger,zemanian} for
applications of the latter to the theory of non random linear
systems. Furthermore, $\bigS_{-1}$ is closed under the Wick
product. More precisely, the following result holds (see
\cite{MR1408433}):
\begin{Tm}[V\aa ge, 1996] \label{Vage}
In the space $\bigS_{-1}= \bigcup_{p} {\mathbf
L}^2(\ell,\mu_{-p})$ it holds that,
\begin{equation}
\label{eqineq} \|f \circ g\|_q \leq A_{q-p} \|f\|_p \|g\|_q,
\end{equation}
(where $\|\cdot\|_p$ denotes the norm of ${\mathbf
L}^2(\ell,\mu_{-p})$) for any $q \geq p+2$, and for any $f \in
\mathbf L_2(\ell, \mu_{-p}),g \in \mathbf L_2(\ell, \mu_{-q})$,
with
\[
A_{q-p}=\left(\sum_{\alpha \in \ell}(2 \mathbb
N)^{-\alpha(q-p)}\right)^{\frac 12}<\infty
\]
\end{Tm}

We note that the finiteness of $A_{q-p}$ was proved by Zhang in
\cite{zhang}. It follows from \eqref{eqineq} that the
multiplication operator
\[
M_f\,\,:\,\,g\mapsto f\circ g
\]
is bounded from the Hilbert space $\mathbf L^2(\ell,\mu_{-q})$
into itself where $f\in \mathbf L^2(\ell,\mu_{-p})$ and $q \geq
p+2$. This also allows us to consider power series. If
$\sum_{n=0}^\infty a_n z^n$ converges in the open disk with
radius $R$, then for any $f \in {\mathbf L}^2(\ell,\mu_{-p})$
with $\| f\|_p<\frac{R}{A_{2}}$,  we obtain
\[
\sum_{n=0}^\infty |a_n|\| f^{\circ n}\|_{p+2} \leq
\sum_{n=0}^\infty |a_n| (A_{q-p}\| f\|_p)^n<\infty,
\]
and hence $\sum_{n=0}^\infty a_n f^{\circ n} \in {\mathbf
L}^2(\ell,\mu_{-(p+2)})$. In this way we are also able to
consider the invertible elements of the algebra $\bigS_{-1}$.
These properties among others, which follows by V\aa ge
inequality, are the key tools for the applications described at the beginning.\\

\subsection{The non-commutative case and an outline of the paper}
In a similar way, the non-commutative white noise space is
defined by the full Fock space
\[
\Gamma(\bigH)=\oplus_{n=0}^\infty\bigH^{\otimes n},
\]
where again, one takes $\bigH_0=\mathbf {\mathbf L}^2(\mathbb
R)$, but other choices of $\bigH_0$ are possible. Denoting by
$\widetilde \ell$ the free (non-commutative) monoid generated by
$\mathbb N$, the space $\widetilde \bigW$ is isometrically
isomorphic to ${\mathbf L}^2(\widetilde \ell, \nu)$, where $\nu$
is now the counting measure (the $\alpha!$ disappeared since we
are no longer in the symmetric case). The non-commutative Wick
product is defined by $(f,g) \mapsto f \otimes g$, and in view of
proposition \ref{ncwick},  $\widetilde \bigW$ is not closed under
it. The counterpart of $\mathcal S_{-1}$ is now of the form
$\bigcup_p\mathbf L^2(\widetilde{\ell},\widetilde{\mu}_{-p})$
where the measures $\widetilde{\mu}_{-p}$ are defined by
\eqref{tildemup}. In the construction of the non-commutative
version of the Kondratiev space of stochastic distributions, an
inequality similar to
the one presented in Theorem \ref{Vage} will be seen to hold.\\

The outline of the paper is as follows: In Section 2 we construct
the non-commutative version of the Kondratiev space, $\widetilde
\bigS_{-1}$. In Section 3, we discuss about second quantization,
and present an inequality which holds in $\widetilde \bigS_{-1}$.
Power series, invertible elements and some other properties
presented in Section 4. In Section 5, we consider $\widetilde
\bigS_{-1}$ as an appropriate framework to stochastic linear
systems.


\section{The white noise space and the Kondratiev space of stochastic
distributions - the non-commutative case} \label{Sec4}
\setcounter{equation}{0}
To define the  non-commutative version of the Gelfand triple
$(\bigS_1,\mathcal W,\bigS_{-1})$, two approaches are possible.
In the first one, we replace the free commutative monoid generated
by $\mathbb N$, namely $\ell$, with the free non-commutative
monoid $\widetilde \ell$ generated by $\mathbb N$. To ease the
notation, we in fact consider a family of (pairwise distinct)
symbols $(z_n)_{n\in\mathbb N}$ indexed by $\mathbb N$, and
consider equivalently the free non-commutative monoid they
generate:
\[
\begin{split}
\widetilde \ell
&= \mathbb N^*\\
&\cong\{ z_{i_1}^{\alpha_{1}}z_{i_2}^{\alpha_{2}} \cdots
z_{i_n}^{\alpha_{n}} : n \in \mathbb N, i_1 \neq i_2 \neq \cdots
\neq i_n
\in \mathbb N, \alpha_1,\dots,\alpha_n \in \mathbb N\}\cup \{1\}\\
&\cong\{ z_{i_1}z_{i_2}\cdots z_{i_m} : m \in \mathbb N, i_1,
\dots,i_n \in \mathbb N \}\cup \{1\}.
\end{split}
\]
We also consider the induced partial order, that is for
$\alpha,\beta \in \widetilde \ell$, we define $\alpha \leq \beta$
if there exists $\gamma
\in \widetilde \ell$ such that $\alpha\gamma=\beta$.\\

For $\alpha=z_{i_1}^{\alpha_{1}}z_{i_2}^{\alpha_{2}} \cdots
z_{i_n}^{\alpha_{n}} \in \widetilde \ell$ (where $i_1 \neq i_2
\neq \cdots \neq i_n$) we define
\[
(2 \mathbb N)^\alpha=\prod_{k=1}^n(2i_k)^{\alpha_k}= \prod_{j \in
\{i_1,\dots,i_n\}}(2j)^{\left( \sum_{k:i_k=j}\alpha_k \right)}.
\]

We define the measures $\widetilde{\nu}(\alpha)=1$ for every
$\alpha\in\widetilde{\ell}$ and for $p\in\mathbb Z$,
\begin{equation}
\label{tildemup} \widetilde \mu_p(\alpha)=(2\mathbb N)^{\alpha p}.
\end{equation}
\begin{Dn}\label{maindef}
We call $\mathbf {\mathbf L}^2(\widetilde \ell,\widetilde\nu)$ the
non-commutative white noise space and we denote it by
$\widetilde{\mathcal W}$. Similarly, $\widetilde \bigS_1 =
\bigcap_{p\in\mathbb N} \mathbf {\mathbf L}^2(\widetilde \ell,
\mu_{p})$ and $\widetilde \bigS_{-1} = \bigcup_{p\in\mathbb N}
\mathbf {\mathbf L}^2(\widetilde \ell, \mu_{-p})$, topologized as
a countably Hilbert space and as its strong dual respectively,
will be called the non-commutative Kondratiev space of stochastic
test functions and the non-commutative Kondratiev space of
stochastic distributions respectively.
\end{Dn}

In the second approach to consider the non-commutative version of
the triple $(\bigS_1,\mathcal W,\bigS_{-1})$ we replace the
symmetric Fock space with the full Fock space. Recall that the
full Fock space over $\bigH$ is the Hilbert space
\[
\Gamma(\bigH)=\oplus_{n=0}^\infty\bigH^{\otimes n}.
\]
Assuming $(e_i)_{i \in I}$ is an orthonormal basis of $\bigH$,
for $\alpha=z_{i_1}^{\alpha_1}z_{i_2}^{\alpha_2} \cdots
z_{i_m}^{\alpha_m} $ (where $i_1 \neq i_2 \neq \cdots \neq i_m
\in I$), such that $|\alpha|=\sum_{j=1}^m \alpha_{j}=n$, we denote
\[
e_\alpha=e_{i_1}^{\otimes \alpha_{1}}\otimes \cdots \otimes
e_{i_m}^{\otimes \alpha_{m}} \in \bigH^{\otimes n}.
\]
$(e_\alpha)$ is clearly an orthonormal basis of $\bigH^{\otimes
n}$ (with respect to the tensor inner product $\langle u_1
\otimes \cdots \otimes u_n , v_1 \otimes \cdots \otimes v_n
\rangle =\prod_{i=1}^n\langle u_i,v_i\rangle $).\\

As in the commutative case we make the choice $\bigH=\mathbf
{\mathbf L}^2(\mathbb R)$ and denote by $(e_n)_{n\in\mathbb N}$ an
orthonormal basis of it (e.g. the Hermite functions). For any $p
\in \mathbb Z$, we denote
\[
\bigH_p=\left\{ \sum_{n=1}^\infty f_n e_n: \sum_{n=1}^\infty
|f_n|^2(2n)^p<\infty\right\} \cong {\mathbf L}^2(\mathbb N,
(2n)^p).
\]
\begin{Rk}{\rm
We note that
\[
 \cdots \subseteq \bigH_2 \subseteq \bigH_1 \subseteq \bigH_0
 \subseteq \bigH_{-1} \subseteq \bigH_{-2}\subseteq \cdots ,
\]
and that $\bigcap_p \bigH_p$ is the Schwartz space of rapidly
decreasing complex smooth functions (in case we indeed choose
$(e_n)$ to be the Hermite functions)  and $\bigcup_p \bigH_p$ is
its dual, namely the Schwartz space of complex tempered
distributions.}
\end{Rk}

\begin{Tm}
It holds that
\[
\widetilde \bigS_1 = \bigcap_{p \in \mathbb N}\Gamma(\bigH_p),
\quad\widetilde \bigW = \Gamma(\bigH_0), \quad \text {and } \quad
\widetilde \bigS_{-1} = \bigcup_{p \in \mathbb
N}\Gamma(\bigH_{-p}).
\]
\end{Tm}
\begin{proof}[Proof]
Clearly $((2n)^{-p/2}e_n)$ is an orthonormal basis of $\bigH_p$.
Hence,
\[
\begin{split}
e^{(p)}_\alpha &=((2i_1)^{-p/2}e_{i_1})^{\circ \alpha_{i_1}}\circ
\cdots \circ
((2i_m)^{-p/2}e_{i_m})^{\circ \alpha_{i_m}}\\
&=\prod_{j=1}^m(2i_j)^{-\alpha_{i_j}p/2}e_\alpha\\
&=(2\mathbb N)^{-\alpha p /2}e_\alpha
\end{split}
\]
is an orthonormal basis of $\Gamma(\bigH_p)$. Thus,
\[
\Gamma(\bigH_{p}) =\left\{ \sum_{\alpha \in \widetilde \ell}
f_\alpha e_\alpha: \sum_{\alpha \in \widetilde \ell}
|f_\alpha|^2(2\mathbb N)^{\alpha p} <\infty \right\},
\]
and so
\[
\begin{split}
\bigcap_{p \in \mathbb N} \Gamma(\bigH_{p}) &=\left\{
\sum_{\alpha \in\widetilde \ell} f_\alpha  e_\alpha: \sum_{\alpha
\in \widetilde\ell} |f_\alpha|^2(2\mathbb N)^{\alpha p}
<\infty \quad \forall p\in\mathbb N \right\} \\
&= \bigcap_p {\mathbf L}^2(\widetilde \ell, \mu_{p})\\
&= \widetilde \bigS_{1},
\end{split}
\]
\[
\Gamma(\bigH_{0}) =\left\{ \sum_{\alpha \in \widetilde\ell}
f_\alpha e_\alpha:\sum_{\alpha \in\widetilde \ell}
|f_\alpha|^2<\infty \right\} = \mathbf {\mathbf L}^2(\widetilde
\ell, \nu)= \widetilde \bigW,
\]
and
\[
\begin{split}
\bigcup_{p \in \mathbb N} \Gamma(\bigH_{-p}) &=\left\{
\sum_{\alpha \in \widetilde\ell} f_\alpha  e_\alpha: \sum_{\alpha
\in \widetilde\ell} |f_\alpha|^2(2\mathbb N)^{-\alpha p}
<\infty \quad \text{for some } p\in\mathbb N \right\} \\
&= \bigcup_p {\mathbf L}^2(\widetilde \ell, \mu_{-p})\\
&=\widetilde \bigS_{-1}.
\end{split}
\]
\end{proof}

As was mentioned in the commutative case, we do not use in this paper any realization of the white noise space.
Similarly to the commutative case, there is an isomorphism between the full Fock space $\Gamma(\bigH_0)$ (i.e. the non-commutative white noise space) and the $\mathbf L^2$-space of the free white noise, namely $\mathbf L^2(\tau)$, where $\tau$ is a free expectation. For more information, we refer to the paper \cite{MR2540072} of M. Bo{\.z}ejko and E. Lytvynov.

\begin{Dn}
The Wick product is defined by $(f,g) \mapsto f \otimes g$
whenever it make sense. In terms of the basis we obtain
\[
f \otimes g = \left(\sum_{\alpha \in \widetilde{\ell}}f_\alpha e_\alpha
\right)\otimes \left(\sum_{\alpha \in \widetilde{\ell}}g_\alpha e_\alpha
\right)=\sum_{\alpha \in \widetilde{\ell}} \left( \sum_{\beta \leq
\alpha}f_\beta g_{\beta^{-1}\alpha} \right) e_\alpha,
\]
where $\beta\le \alpha$ means that there exists there exists (a unique) $\gamma\in\widetilde{\ell}$  such that
$\alpha=\beta\gamma$, and $\beta^{-1}\alpha$ stands for $\gamma$.
\end{Dn}
Thus, the Wick product is the convolution of functions over the
monoid $\widetilde \ell$.

\begin{Pn} \label{ncwick}
$\widetilde \bigW$ is not closed under the Wick product.
\end{Pn}
\begin{proof}[Proof]
Let $\iota: \ell^2(\mathbb N) \to \widetilde W$ be the embedding
defined by
\[
\langle \iota(f), e_\alpha \rangle =
 \begin{cases}
   f_n & \text{if } \alpha=z_1^n\\
   0 & \text{otherwise}
  \end{cases}
\]
(where $f=(f_n) \in \ell^2(\mathbb N)$), and let $f,g \in
\ell^2(\mathbb N)$ such that $\|f*g\|=\infty$, where $*$ denotes
the standard convolutions of two elements in $\ell^2(\mathbb N)$.
Then,
\[
\| \iota(f) \otimes \iota(g) \|=\|f*g\|=\infty.
\]
\end{proof}
\begin{Rk}\label{notclosed}{\rm
The reason why the commutative white noise space is not closed
under the symmetric Wick product is similar. We can simply define
$\eta: \ell^2(\mathbb N) \to W$ by
\[
\langle \eta(f), e_\alpha \rangle =
 \begin{cases}
   f_n/{\sqrt{n!}} & \text{if } \alpha=(n,0,0,\ldots)\\
   0 & \text{otherwise}
  \end{cases}
\]
(where $f=(f_n) \in \ell^2(\mathbb N)$). Thus, for non-negative
sequences $f,g \in \ell^2(\mathbb N)$ such that $\|f*g\|=\infty$,
\[
\begin{split}
\| \eta(f) \otimes \eta(g) \|^2
&= \sum_n \left( \sum_{k=1}^n \frac{1}{\sqrt{k!(n-k)!}} f_k g_{n-k} \right)^2 n!\\
&\geq \sum_n \left( \sum_{k=1}^n f_k g_{n-k} \right)^2 \\
&= \|f*g\|^2\\
&=\infty.
\end{split}
\]}
\end{Rk}

Similar to the commutative case, it will be shown in the sequel
that $\widetilde \bigS_{-1}$ is closed under the Wick product, and
moreover it satisfies an inequality similar to the one that was
presented in Theorem \ref{Vage}.

\section{Second quantization and an inequality of tensor product}\label{Sec5}
\setcounter{equation}{0} Let $\bigK_0$ be a separable Hilbert
space, and let $(e_n)_{n\in\mathbb N}$ be an orthonormal basis of
$\bigK_0$. Furthermore, let $(a_n)_{n\in\mathbb N} $ be a
sequence of real numbers greater than or equal to 1. For any $p \in \mathbb Z$, we denote
\[
\bigK_p=\left\{ \sum_{n=1}^\infty f_n e_n: \sum_{n=1}^\infty
|f_n|^2 a_n^p< \infty\right\} \cong {\mathbf L}^2(\mathbb N,
a_n^p).
\]
We note that
\[
\cdots \subseteq \bigK_2 \subseteq \bigK_1 \subseteq \bigK_0
\subseteq \bigK_{-1} \subseteq \bigK_{-2}\subseteq \cdots,
\]
where the embedding $T_{q,p}:\bigK_q \hookrightarrow \bigK_p$
satisfies
\[
\|T_{q,p} a_n^{-q/2}e_n\|_p =a_n^{-(q-p)/2}\| a_n^{-p/2} e_n\|_q,
\]
and hence
\[
\|T_{q,p}\|_{HS}=\sqrt{\sum_{n \in \mathbb N} a_n^{-(q-p)}}.
\]
The dual of a Fr\'echet space is nuclear if and only if the
initial space is nuclear. Thus, $\bigcup_{p \in \mathbb N}
\bigK_{-p}$ is nuclear if and only if $\bigcap_{p \in \mathbb N}
\bigK_p$ is nuclear. This is turn will hold if and only if for any
$p$ there is some $q>p$ such that $\|T_{q,p}\|_{HS}<\infty$, that
is, if and only if there exists some $d>0$ such that $\sum_{n \in
\mathbb N} a_n^{-d}$ converges. We note that in this case, $d$
can be chosen so that
\[
\sum_{n \in \mathbb N} a_n^{-d}<1.
\]
We call the smallest integer $d$ which satisfy this inequality
the index of $\bigcup_{p \in \mathbb N} \bigK_{-p}$. In this
section we show that if $\bigcup_{p \in \mathbb N} \bigK_{-p}$ is
nuclear of index $d$, then  $\bigcup_{p \in \mathbb
N}\Gamma(\bigK_{-p})$ has the property that
\[
\|f \otimes g\|_q \leq \|\Gamma(T_{q,p})\|_{HS}\|f\|_p\|g\|_q \text{ and }
\|g \otimes f\|_q \leq \|\Gamma(T_{q,p})\|_{HS}\|f\|_p\|g\|_q
\]
for all $q \geq p+d$,
where $\|\cdot\|_p$ is the norm associated to
$\Gamma(\bigK_{-p})$, and $\|\Gamma(T_{q,p})\|_{HS}$ is finite.
The case $a_n=2n$ (and hence $d=2$) corresponds to the
non-commutative Kondratiev space, and is discussed in the next
section.\\

\begin{Dn}\label{2ndQuantization}
Let $T:\bigH_1 \to \bigH_2$ be a bounded linear operator between
two separable Hilbert spaces. Then $T^{\otimes n}:\bigH_1^{\otimes n} \to
\bigH_2^{\otimes n}$, defined by
\[
T^{\otimes n}(u_1 \otimes \cdots \otimes u_n)=Tu_1 \otimes \cdots
\otimes Tu_n,
\]
is a bounded linear operator between $\bigH_1^{\otimes n}$ and
$\bigH_2^{\otimes n}$.  When $T$ is a contraction, it induces a
bounded linear operator $\Gamma(\bigH_1) \to \Gamma(\bigH_2)$,
denoted by $\Gamma(T)$, and called the second quantization of $T$.
\end{Dn}

Let $(\lambda_n)$ be a sequence of non-negative numbers. For
$\alpha=z_{i_1}^{\alpha_{1}}z_{i_2}^{\alpha_{2}} \cdots
z_{i_n}^{\alpha_{n}} \in \widetilde \ell$ (where $i_1 \neq i_2
\neq \cdots \neq i_n$) we denote
\[
\lambda_{\mathbb
N}^\alpha=\prod_{k=1}^n\lambda_{i_k}^{\alpha_k}=\prod_{j \in
\{i_1,\dots,i_n\}}\lambda_j^{\left( \sum_{k:i_k=j}\alpha_k
\right)}.
\]
We recall that if $T:\bigH_1 \to \bigH_2$ is a compact operator between two separable Hilbert spaces,
then
\[
Tf=\sum_{n=1}^\infty\lambda_n \langle f, e_n \rangle h_n
\]
where $(e_n)_{n\in\mathbb N}$ and $(h_n)_{n\in\mathbb N}$ are
orthonormal basis of $\bigH_1$ and $\bigH_2$ respectively and
where $(\lambda_n)$ is a non-negative sequence converging to zero.
Conversely, any such a decomposition defines a compact operator
$\bigH_1 \to \bigH_2$ (see for instance \cite{MR58:12429a}).

\begin{Tm}
\label{thm:HS} Let $T:\bigH_1 \to \bigH_2$ be a compact contraction operator between two separable Hilbert spaces
with
\[
Tf=\sum_{n=1}^\infty\lambda_n \langle f, e_n \rangle h_n
\]
where $(e_n)_{n\in\mathbb N}$ and $(h_n)_{n\in\mathbb N}$ are
orthonormal basis of $\bigH_1$ and $\bigH_2$ respectively and
where $(\lambda_n)$ is a non-negative sequence converging to
zero. Let $\Gamma(T)$ be its second quantization as in Definition \ref{2ndQuantization}.
Then,
\begin{enumerate}[(a)]
\item It holds that
\[
\Gamma(T)f=\sum_{\alpha \in \widetilde \ell} \lambda_{\mathbb
N}^{\alpha} \langle f, e_\alpha \rangle h_\alpha,
\]
where $(e_\alpha)_{\alpha\in\widetilde{\ell}}$ and
$(h_\alpha)_{\alpha\in\widetilde{\ell}}$ are orthonormal basis of
$\Gamma(\bigH_1)$ and $\Gamma(\bigH_2)$ respectively.
\item If furthermore $T$ is an Hilbert-Schmidt operator, i.e.
$(\lambda_n) \in \ell^2(\mathbb N)$, then
\[
\|\Gamma(T)\|_{HS}^2=\sum_{n=0}^\infty\|T\|_{HS}^{2n}.
\]
In particular, $\Gamma(T)$ is a Hilbert-Schmidt operator if and
only if $T$ is a Hilbert-Schmidt operator with $\|T\|_{HS}<1$ and
in this case we obtain
\[
\|\Gamma(T)\|_{HS}=\frac 1{\sqrt{1-\|T\|_{HS}^2}}
\]
\end{enumerate}
\end{Tm}

\begin{proof}[Proof]
For any $\alpha \in \widetilde \ell$  let
$e_\alpha=e_{i_1}^{\otimes \alpha_{1}}\otimes \cdots \otimes
e_{i_m}^{\otimes \alpha_{m}}$ and $h_\alpha=h_{i_1}^{\otimes
\alpha_{1}}\otimes \cdots \otimes h_{i_m}^{\otimes \alpha_{m}}$.
Then, $(e_\alpha)_{\alpha\in\widetilde{\ell}}$ and
$(h_\alpha)_{\alpha\in\widetilde{\ell}}$ are orthonormal basis of
$\bigH_1$ and $\bigH_2$ respectively.
\begin{enumerate}[(a)]
\item We have that
\[
\begin{split}
\Gamma(T)e_\alpha &=(Te_{i_1})^{\otimes \alpha_{1}}\otimes \cdots
\otimes
(Te_{i_m})^{\otimes \alpha_{m}}\\
&=(\lambda_{i_1}h_{i_1})^{\otimes \alpha_{1}}\otimes \cdots
\otimes (\lambda_{i_m}h_{i_m})^{\otimes \alpha_{m}}\\
&=\lambda_{\mathbb N}^\alpha h_\alpha.
\end{split}
\]
Thus, by the linearity and continuity of $\Gamma(T)$,
\[
\Gamma(T)f=\sum_{\alpha \in \widetilde \ell} \lambda_{\mathbb
N}^{\alpha} \langle f, e_\alpha \rangle h_\alpha.
\]

\item
We have that
\begin{align*}
\|\Gamma(T)\|_{HS}^2
&=\sum_{\alpha \in \widetilde \ell}\|\Gamma(T)e_\alpha\|^2\\
\displaybreak[1]&=\sum_{n=0}^\infty\sum_{\alpha \in \widetilde
\ell,|\alpha|=n}
\|T^{\otimes n}e_\alpha\|^2\\
\displaybreak[1]&=\sum_{n=0}^\infty\sum_{\alpha \in \widetilde
\ell,|\alpha|=n}
\prod_{i=1}^\infty \|Te_i\|^{2\alpha_i}\\
\displaybreak[1]&=\sum_{n=0}^\infty\sum_{\alpha \in \widetilde{\ell}
,|\alpha|=n}\frac{n!}{\alpha!}\prod_{i=1}^\infty
\|Te_i\|^{2\alpha_i}.
\end{align*}
Considering an experiment with $\mathbb N$ results, where the
probability of the result $i$ is $p_i=\|T\|_{HS}^{-2}\|Te_i\|^2$
(and so $\sum p_i=1$), the probability that repeating the
experiment $n$ times yields that the result $i$ occurs $\alpha_i$
times for any $i$ is
\[
\frac{n!}{\alpha!}\prod_{i=1}^\infty p_i^{\alpha_i}
=\|T\|_{HS}^{-2n}\frac{n!}{\alpha!}\prod_{i=1}^\infty
\|Te_i\|^{2\alpha_i}.
\]
Thus,
\[
\sum_{\alpha \in \widetilde{\ell}
,|\alpha|=n}\frac{n!}{\alpha!}\prod_{i=1}^\infty
\|Te_i\|^{2\alpha_i}=\|T\|_{HS}^{2n},
\]
and we obtain the requested result.
\end{enumerate}
\end{proof}

\begin{Tm}\label{Vage2}
If  $\bigcup_{p \in \mathbb N} \bigK_{-p}$ is nuclear of index
$d$, then  $\bigcup_{p \in \mathbb N}\Gamma(\bigK_{-p})$ is
nuclear and has the property that
\[
\|f \otimes g\|_q \leq \|\Gamma(T_{q,p})\|_{HS}\|f\|_p\|g\|_q \text{ and }
\|g \otimes f\|_q \leq \|\Gamma(T_{q,p})\|_{HS}\|f\|_p\|g\|_q
\]
for all $q\geq p+d$, where $\|\cdot\|_p$ is the norm associated to
$\Gamma(\bigK_{-p})$, and where
\[
\|T_{q,p}\|_{HS}=\sum_{\alpha \in \widetilde \ell}a_{\mathbb
N}^{-\alpha (q-p)}=\frac 1{\sqrt{1-\sum_{n \in \mathbb
N}a_n^{-(q-p)}}}.
\]
\end{Tm}
\begin{proof}
Denoting $b_{\alpha}=a^\alpha_{\mathbb N}$, we have that
\[
\Gamma(\bigK_{-p})=\left\{(f_\alpha)_{\alpha \in \widetilde \ell}:
\sum_{\alpha \in \widetilde \ell}
|f_\alpha|^2b_\alpha^{-p}<\infty\right\}.
\]
Since $\bigcup_{p \in \mathbb N} \bigK_{-p}$ is nuclear of index
$d$,
\[
\|T_{q,p}\|^2=\sum_{n \in \mathbb N}a_n^{-(q-p)}<1 \quad\text{
for any }q\geq p+d
\]

In view of Theorem \ref{thm:HS}, $\Gamma(T_{q,p})$ is
Hilbert-Schmidt and
\[
{\sum_{\alpha \in \widetilde \ell}
b_{\alpha}^{-(q-p)}}={\sum_{\alpha \in \widetilde \ell} a_{\mathbb
N}^{-\alpha(q-p)}}=\|\Gamma(T_{q,p})\|_{HS}^2=\frac
1{{1-\|T_{q,p}\|_{HS}^2}}<\infty.
\]
Since for any $\alpha=z_{i_1}^{\alpha_{1}}z_{i_2}^{\alpha_{2}}
\cdots z_{i_n}^{\alpha_{n}} \in \widetilde \ell$ and
$\beta=z_{j_1}^{\beta_{1}}z_{j_2}^{\beta_{2}} \cdots
z_{i_m}^{\beta_{m}} \in \widetilde \ell$ it holds that
\[
b_{\alpha}b_{\beta} =a^\alpha_{\mathbb N}a^\beta_{\mathbb N}
=\prod_{k=1}^n a_{i_k}^{\alpha_k}\cdot \prod_{l=1}^m
a_{i_l}^{\beta_l} =a^{\alpha \beta}_{\mathbb N} =b_{\alpha\beta},
\]
for any $f \in \Gamma(\bigH_{-p})$ and $g \in \Gamma(\bigH_{-q})$
we obtain
\begin{align*}
\|f \otimes g\|_{q} ^2 & =\sum_{\gamma \in\widetilde
\ell}\left|\sum_{\alpha \leq \gamma} f_\alpha
g_{\alpha^{-1}\gamma}
b_\gamma^{-q/2}\right|^2\\
\displaybreak[1]& \leq \sum_{\gamma \in \widetilde
\ell}\left(\sum_{\alpha \leq \gamma}| f_\alpha| b_{\alpha}^{-q/2}
|g_{\alpha^{-1}\gamma}|
b_{\alpha^{-1}\gamma}^{-q/2}\right)^2\\
\displaybreak[1]& = \sum_{\gamma \in \widetilde
\ell}\left(\sum_{\alpha,\alpha' \leq \gamma}|
f_\alpha|b_{\alpha}^{-q/2}|f_{\alpha'}| b_{\alpha'}^{-q/2}
|g_{\alpha^{-1}\gamma}|b_{\alpha^{-1}\gamma}^{-q/2}|g_{(\alpha')^{-1}\gamma}
|b_{(\alpha')^{-1}\gamma}^{-q/2}\right)\\
\displaybreak[1]& \leq \sum_{\alpha,\alpha' \in \widetilde
\ell}\left(| f_\alpha|b_{\alpha}^{-q/2} |f_{\alpha'}|
b_{\alpha'}^{-q/2} \sum_{\gamma \geq \alpha,\alpha'}
|g_{\alpha^{-1}\gamma}|b_{\alpha^{-1}\gamma}^{-q/2}|g_{(\alpha')^{-1}\gamma}
|b_{(\alpha')^{-1}\gamma}^{-q/2} \right)\\
\displaybreak[1]& \leq \left(\sum_{\beta \in \widetilde
\ell}|f_\beta| b_{\beta}^{-p/2} \right)^2 \left(\sum_{\beta \in
\widetilde \ell}|g_{\beta}|^2 b_{\beta}^{-q}\right)^{\frac 12}
 \left(\sum_{\beta \in \widetilde \ell}|g_{\beta}|^2 b_{\beta}^{-q}\right)^{\frac 12}\\
\displaybreak[1]&\leq \left(\sum_{\beta \in \widetilde \ell}
b_{\beta}^{-(q-p)}\right) \left(\sum_{\beta \in \widetilde \ell}
|f_\beta|^2 b_{\beta}^{-p}\right) \left(\sum_{\beta
 \in \widetilde \ell}|g_{\beta}|^2 b_{\beta}^{-q}\right)\\
\displaybreak[1]&=\|
\Gamma(T_{q,p})\|_{HS}^2\|f\|_{p}^2\|g\|_{q}^2.
\end{align*}
The second inequality is obtained in the same manner since
\[
(f\otimes g)_\gamma=
\sum_{\substack{\alpha,\beta\in\widetilde{\ell}\\ \alpha\beta=\gamma}} f_\alpha
g_{\beta}=\sum_{\substack{\alpha,\beta\in\widetilde{\ell}\\ \beta\alpha=\gamma}}
 f_{\beta}
g_{\alpha}.
\]
\end{proof}

\section{The algebra of the non-commutative Kondratiev space of
stochastic distributions}\label{Sec6} \setcounter{equation}{0}
We now specialize the results of the preceding section to
$a_n=2n$, and denote by $\bigH_p$ the corresponding spaces:
\[
\bigH_{p}=\left\{ \sum_{n=1}^\infty f_n e_n: \sum_{n=1}^\infty
|f_n|^2(2n)^{p}<\infty\right\} \cong \mathbf {\mathbf
L}^2(\mathbb N, (2n)^{p}),
\]
Denoting by $T_{q,p}$ the embedding $\bigH_q \hookrightarrow
\bigH_p$, it holds that
\[
\|T_{q,p}\|_{HS}^2=\sum_{n \in \mathbb N}
(2n)^{-(q-p)}=2^{-(q-p)} \zeta(q-p),
\]
where $\zeta$ denotes Riemann's zeta function. Since for any $s
\geq 2$, $\zeta(s) < 2^{s}$, for any $q \geq p+2$,
$\|T_{q,p}\|_{HS}<1$. In view of Theorems \ref{thm:HS} and
\ref{Vage2} we obtain the following result:\mbox{}\\

\begin{Tm}\label{mainthm}
\begin{enumerate}[(a)]
\item
The non-commutative Kondratiev spaces $\widetilde \bigS_1$ and
$\widetilde \bigS_{-1}$ are both nuclear spaces.
\item For any $q \geq p+2$,
\[
B_{q-p}^2={\sum_{\alpha \in \widetilde \ell} (2\mathbb
N)^{-\alpha(q-p)}}= \frac 1{{1-2^{-(q-p)} \zeta(q-p)}},
\]
where $B_{q-p}=\| \Gamma(T_{q,p})\|_{HS}$.
\item For any $q \geq p+2$ and for any $f \in \Gamma(\bigH_{-p})$  and
$g \in \Gamma(\bigH_{-q})$
\begin{equation}
\label{vage2} \|f \otimes g\|_q \leq B_{q-p} \|f\|_p \|g\|_q \quad \text{ and }  \quad \|g \otimes f\|_q \leq B_{q-p} \|f\|_p \|g\|_q
\end{equation}
where $\|\cdot\|_p$ is the norm associated to
$\Gamma(\bigH_{-p})$.
\end{enumerate}
\end{Tm}

We now show that the non-commutative Wick product is continuous.
We first need the following proposition.

\begin{Pn}
\label{pwcont}
Let $f \in \widetilde \bigS_{-1}$.
Then the linear mappings $L_a:x \mapsto ax$, $R_a:x \mapsto xa$ are continuous.
\end{Pn}
\begin{proof}[Proof]
Suppose that $f \in \Gamma(\bigH_{-p})$, and
let $L_a|_{\Gamma(\bigH_{-r})} : \Gamma(\bigH_{-r}) \to \widetilde \bigS_{-1}$ be the restriction of the map $L_a$ to $\Gamma(\bigH_{-r})$.
If $B$ is a bounded set of $\Gamma(\bigH_{-r})$ then in particular we may choose $q \geq p+2$ such that $q \geq r$, so
$B \subseteq \{ g \in \Gamma(\bigH_{-q}): \|g\|_q <\lambda \}$.
Thus, for any $g \in B$
\[
\|L_a|_{\Gamma(\bigH_{-q})} (g) \|_q \leq B_{q-p}\lambda \|g\|_q.
\]
Hence, $L_a|_{\Gamma(\bigH_{-q})}(B)$ is bounded in $\Gamma(\bigH_{-q})$ and hence in $\widetilde \bigS_{-1}$. Thus, for any $r$, $L_a|_{\Gamma(\bigH_{-r})} :\Gamma(\bigH_{-r}) \to \widetilde \bigS_{-1}$ is bounded and hence continuous.
Since $\widetilde \bigS_{-1}=\bigcup_{p \in \mathbb N} \Gamma(\bigH_{-p})$ is a strong dual of the reflexive Fr\'echet space $\widetilde \bigS_{1}=\bigcap_{p \in \mathbb N} \Gamma(\bigH_p)$, it is the inductive limit of the Hilbert spaces $\Gamma(\bigH_{-p})$ (see \cite[IV.23]{BourbakiTVS}).
So by the universal property of inductive limits, $L_a$ is continuous. The proof for $R_a$ is similar.
\end{proof}

\begin{Tm}
The Wick product is a continuous function $\widetilde \bigS_{-1}
\times \widetilde \bigS_{-1} \to \widetilde \bigS_{-1}$ in the
strong topology. Hence $(\widetilde \bigS_{-1},+,\otimes)$ is a
topological $\mathbb C$-algebra.
\end{Tm}

This follows immediately from Proposition \ref{pwcont} together with the following theorem, proved in
\cite[IV.26]{BourbakiTVS}.
\begin{Tm}
Let $E_1$ and $E_2$ be two reflexive Fr\'echet spaces, and let $G$ a locally convex Hausdorff space. For $i=1,2$, let $F_i$ be the strong dual of $E_i$. Then every separately continuous bilinear mapping $u:F_1 \times F_2 \to G$ is continuous.
\end{Tm}

As a matter of fact, the topology of the space $\widetilde \bigS_{-1}$ itself is hardly used, and most of the applications only its ``local topology", i.e. the topology of the Hilbert spaces $\Gamma(\bigH_{-p})$.
Nonetheless, we give here, as a remark, a brief discussion about its topology and about the relations of this topology to the topologies of the Hilbert spaces $\Gamma(\bigH_{-p})$.
\begin{Rk}
{\rm $\widetilde \bigS_{-1}$ carries out {\it a priori} two
natural topologies. The first is its topology as a strong dual of
Fr\'echet space (namely, $\widetilde \bigS_{1}$). This topology
was in our mind during our discussion up to now (see Definition
\ref{maindef}). Two of the main properties of this topology is
that any bounded set of $\widetilde \bigS_{-1}$, is bounded in
some Hilbert space $\Gamma(\bigH_{-p})$, and that if the
Fr\'echet space is nuclear (as in our case), then so is its
strong dual (see Theorem \ref{mainthm}(a) and (see
\cite[IV.21-26, \S 3]{BourbakiTVS} and \cite[ \S 5]{GS2_english}
for references on this fact).\smallskip

 The second topology is
its topology as an inductive limit of the locally convex spaces
(which are actually Hilbert spaces) $\Gamma(\bigH_{-p})$, i.e.
the finest locally convex topology such that the embeddings
$\Gamma(\bigH_{-p}) \hookrightarrow \widetilde \bigS_{-1}$ are
continuous.
 There are two main properties of this topology which are worth mentioning.
The first is that it satisfies the universal property of an
inductive limit, i.e. any linear map from an inductive limit of a
family of locally convex spaces to another locally convex space
is continuous if and only if the restriction of the map to any of
members of the family is continuous (see
\cite[II.29]{BourbakiTVS}). The second property is that in case
the inductive limit is of Banach spaces (recall that in our case
they are Hilbert spaces), then the inductive limit is
bornological (see  \cite[III.11-13, \S 2]{BourbakiTVS}) and
barreled (see  \cite[III.24-25, \S 4]{BourbakiTVS}). In our case,
where the ``building block" spaces $\Gamma(\bigH_{-p})$ are
Hilbert spaces (actually, reflexive Banach spaces is enough),
these two topologies coincide (see the proof of \cite[IV.23,
Proposition 4.]{BourbakiTVS}). Furthermore, a nice property
holds: since the embeddings of $\Gamma(\bigH_{-p})$ in
$\Gamma(\bigH_{-q})$ for any $q \geq p+2$ are compact (see
Theorem \ref{mainthm}(b), where it is stated that they are
nuclear, so in particular compact), the topology of  $\widetilde
\bigS_{-1}$ is the finest topology (rather than the finest
locally convex topology) such that the embeddings
$\Gamma(\bigH_{-p}) \hookrightarrow \widetilde \bigS_{-1}$ are
continuous (see \cite[III.6, Lemma 1.]{BourbakiTVS}). }
\end{Rk}

Note that the topological $\mathbb C$-algebra $(\widetilde \bigS_{-1},+,\otimes)$ is
unital, where the unit element is $e_0=\mathbf 1$ which is also the vacuum
vector of $\widetilde \bigW$ embedded in $\widetilde \bigS_{-1}$.
\begin{Dn}
Let $f =\sum_{\alpha \in \widetilde \ell}f_\alpha e_\alpha \in
\widetilde\bigS_{-1}$. Then, $f_0 \in \mathbb{C}$ is called the
generalized expectation of $f$ and is denoted by $E[f]$.
\end {Dn}
From this definition we have
\[
E[f \otimes g] = E[f]E[g]  \quad{\rm and} \quad E[{\mathbf 1}]=1 \quad \quad \forall f,g \in \bigS_{-1}.
\]
Thus, $E: \widetilde \bigS_{-1} \to \mathbb{C}$ is a unital
algebra homomorphism.
In the sequel, we will see it is
the only homomorphism with this property (see Proposition \ref{prop: homo}).
Note also that for any $p \in \mathbb F$, $|E(f)| =|f_0| \leq \|f\|_p$. Since as the strong dual of the reflexive Fr\'echet space $\widetilde \bigS_{1}=\bigcap_{p \in \mathbb N} \Gamma(\bigH_p)$ ,$\widetilde \bigS_{-1}$ is the inductive limit of the Hilbert spaces $\Gamma(\bigH_{-p})$,
by the universal property of inductive limits, $E$ is continuous.

\begin{Pn} \label{Pn:limF}
For any $f  \in \widetilde
\bigS_{-1}$ such that $E[f]=0$, it holds that $\lim_{q \to
\infty}\|f\|_q=0$.
\end{Pn}

\begin{proof}[\indent Proof]
Let $f =\sum_{\alpha \in \widetilde \ell} f_\alpha e_\alpha \in
\Gamma(\bigH_{-p})$ with $f_0=0$. Then for all $\alpha \in
\widetilde \ell$ we have
\[
\lim_{q \to \infty}|f_\alpha|^2(2 \mathbb N)^{-q\alpha}=0,
\]
and for all $q>p$,
\[
|f_\alpha|^2(2\mathbb N)^{-q\alpha} \leq |f_\alpha|^2(2\mathbb
N)^{-p\alpha},
\]
where $\sum_{\alpha \in \widetilde
\ell}|f_\alpha|^2a_\alpha^{-p}=\|f\|_p^2<\infty.$ Thus, the
dominated convergence theorem implies
\[
\lim_{q \to \infty}\|f\|_q^2 =\lim_{q \to \infty}
\sum_{\alpha\in\widetilde \ell} |f_\alpha|^2 (2\mathbb
N)^{-q\alpha} =\sum_{\alpha\in\widetilde \ell}\lim_{q \to \infty}
|f_\alpha|^2 (2\mathbb N)^{-q\alpha} =0.
\]
\end{proof}

\begin{Pn} \label{Pn:power}
Let $f$ be in $\Gamma(\bigH_{-p})$. Then
\[
f^{\otimes n} \in \Gamma(\bigH_{-(p+2)})\quad \forall n \in
\mathbb{N}.
\]
Moreover,
\[
\|f^{\otimes n}\|_{p+2} \leq B_2^n \|f\|_{p}^n.
\]
\end{Pn}

\begin{proof}[\indent Proof]
Obviously, $f^0=1 \in \Gamma(\bigH_{-(p+2)})$, and
$\|f^{ 0}\|_{p+2} = A(2)^0  \|f\|_{p}^0$.\\
By induction,
\[
\begin{split}
\|f^{\otimes (n+1)}\|_{p+2}
&=\|f \otimes f^{\otimes n}\|_{p+2}\\
&\leq B_2 \|f\|_{p} \|f^{\otimes n}\|_{p+2} \\
&\leq B_2^n \|f\|_{p}^{n+1} <\infty
\end{split}
\]
\end{proof}

More generally, given a polynomial $p(z)=\sum_{n=0}^N p_n z^n$
($p_n \in \mathbb{C}$), we define its Wick version $ p: \widetilde
\bigS_{-1} \rightarrow\widetilde \bigS_{-1}$ by
\[
p(f)=\sum_{n=0}^N p_n f^{\otimes n}
\]
By Proposition \ref{Pn:power}, we have that $p(f) \in \widetilde
\bigS_{-1}$ for $f \in \widetilde \bigS_{-1}$. The following
proposition considers the case of power series.

\begin{Pn}
\label{Pn:series} Let $\phi(z)=\sum_{n=0}^\infty \phi_n
z^n$ be a power series (with complex coefficients) which
converges absolutely in the open disk with radius $R$. Then for
any $f \in \widetilde \bigS_{-1}$ such that
$|E[f]|<\frac{R}{B_2}$  it holds that
\[
\phi(f)=\sum_{n =0}^\infty \phi_n f^{\otimes n} \in
\widetilde \bigS_{-1}.
\]
\end{Pn}

\begin{proof}[\indent Proof]
Applying Proposition \ref{Pn:limF}, there exists $q$ such that
\[
\|f-E(f)\|_q<\frac{R}{B_2}-|E[f]|.
\]
Therefore,
\[
\|f\|_q \leq \|f-E(f)\|_q+|E(f)|<\frac{R}{B_2}.
\]
By Proposition \ref{Pn:power}, for all $p \geq q+2$,
\[
\begin{split}
\sum_{n =0}^\infty |\phi_n|\| f^{\otimes n}\|_{p}
&\leq \sum_{n =0}^\infty |\phi_n|B_2 ^n \|f\|_{q}^n\\
&= \sum_{n =0}^\infty |\phi_n|(B_2 \|f\|_{q})^n\\
&<\infty.
\end{split}
\]
Since $\Gamma(\bigH_{-p})$ is a Hilbert space, $\phi(f)=\sum_{n =0}^\infty \phi_n f^{\otimes n} \in \Gamma(\bigH_{-p})$.
Thus, $\phi(f) \in \widetilde \bigS_{-1}$.
\end{proof}

\begin{Pn} \label{prop: inv}
An element $f \in \Gamma(\bigH_{-p})$ is invertible  if and only
if $E[f]$ is invertible.
\end{Pn}

\begin{proof}[\indent Proof]
If $E[f] \neq 0$, we can assume that $E[f] =1$. By Proposition
\ref{Pn:series} we have that $\sum_{n =0}^\infty
(1-f)^{\otimes n} \in \widetilde \bigS_{-1}$. Furthermore,
\[
f\otimes \left(\sum_{n =0}^\infty (1-f)^{\otimes n}\right)=1.
\]
Conversely, assume $f$ invertible. Then there exists $f^{-1} \in
\widetilde \bigS_{-1}$ such that $f \otimes  f^{-1} =1$. Hence,
$E[f]E[f^{-1}]=E[f\otimes  f^{-1}]=1$.
\end{proof}

\begin{Pn}\label{prop: homo}
The following properties hold:
\begin{enumerate}[(a)]
\item The set of all invertible elements in $\widetilde \bigS_{-1}$, denoted by $GL(\widetilde \bigS_{-1})$, is open.
\item The spectrum of $f \in \widetilde \bigS_{-1}$, $\sigma(f)=
\{\lambda \in \mathbb C : f-\lambda  \text{ is not invertible }\}$
 is the singleton $\{E[f]\}$.
\item $E$ is the only homomorphism $\widetilde \bigS_{-1} \to \mathbb{C}$ which is unital.
\end{enumerate}
\end{Pn}
\begin{proof}[\indent Proof]
$\quad$\\
\begin{enumerate}[(a)]
\item
By Proposition \ref{prop: inv}, we have  that $\{f \in \widetilde
\bigS_{-1} : E[f] \neq 0\}$ is the set of all invertible elements
in $\widetilde \bigS_{-1}$. In other words, $GL(\widetilde
\bigS_{-1})$ is the inverse image of $GL(\mathbb C)$ under the generalized expectation $E$. In particular, since $E$ is
continuous, $GL(\widetilde \bigS_{-1})$ is open.
\item
Clearly, $f-\lambda { 1}$ does not have an inverse if and only if
$\lambda=E(f)$.
\item
Let $\varphi:{\widetilde \bigS_{-1}} \to \mathbb C$ be a
homomorphism mapping $\mathbf 1$ to $1$,
and let $f \in {\widetilde \bigS_{-1}}$. Since $\varphi
\left(f-\varphi(f) \right)=0$, $\varphi(f) \in \sigma(f)$, that
is $\varphi(f)=E[f]$.
\end{enumerate}
\end{proof}

\section{Applications to non-commutative linear systems}
\label{Sec7} \setcounter{equation}{0}
We refer to \cite{MR839186,MR1071708,MR0255260,SSR} for general
information on the theory of linear systems, including over
commutative rings, and to the papers \cite{MR1168501,MR0452844}
for more information on linear system on non-commutative rings,
and in particular for the notions of controllable and observable
pairs. In the present setting an input-output system will be a
map of the form now an input-output relation of the form
\begin{equation}
\label{acap1}
                y_n=\sum_{m=0}^{n}{h_m \otimes u_{n-m}},\quad
                n\in\mathbb N_0,
                \end{equation}
where the input sequence $(u_n)_{n\in\mathbb N_0}$, the impulse
response $(h_n)_{n\in\mathbb N_0}$ belong to
$\widetilde{S}_{-1}^{q\times 1}$ and $\widetilde{S}_{-1}^{p\times
q}$ respectively. Then, the output sequence belongs to
$\widetilde{S}_{-1}^{p\times 1}$. When the impulse response
$(h_n)$ or the input sequence $(u_n)$ are not random, the Wick
product reduces to the pointwise product of complex numbers, and
we recover classical convolution systems. The transfer function
of the system \eqref{acap1} is (the possibly divergent) series
defined by
\[
\mathscr H(z)=\sum_{n=0}^\infty h_nz^n,
\]
where $z$ is a complex variable. The realization problem in this
setting is to find, when possible, realization of $\mathscr H$ in
the form
\begin{equation}
\label{real2012} \mathscr H(z)=D+zC\otimes (I-zA)^{-1}B,
\end{equation}
where $A,B,C$ and $D$ are matrices of appropriate entries and
with entries in $\widetilde{S}_{-1}$, and
\[
(I-z A)^{-1}=\sum_{k=0}^\infty z^k A^{\otimes k}.
\]
The series converges in a neighborhood of the origin thanks to
Proposition \ref{Pn:series}.\\

The results presented in \cite{aa_goh,alp,MR2610579} for the case
of the commutative Kondratiev space $\mathcal S_{-1}$ of
stochastic distributions still hold for the non-commutative case
because of the underlying structure and in particular of
inequality \eqref{vage2}. We will present here one representative
result, see Theorem \ref{tm:obs}. Note that the arguments in
\cite{aa_goh,alp,MR2610579} are in the setting of power series
(because one considers there the Hermite transform of the
Kondratiev space rather than the Kondratiev space itself), and
make use of derivatives. For the general case, when no power
series are available, we need to introduce  and prove the
continuity, of the operators $D_m$, $m=1,2,\ldots$ defined by
\[
D_m({z_{i_1}^{\alpha_{1}}z_{i_2}^{\alpha_{2}} \cdots
z_{i_n}^{\alpha_{n}}})=
   \sum_{\{j:i_j=m,\alpha_j>0\}}\alpha_j {z_{i_1}^{\alpha_{1}}z_{i_2}^{\alpha_{2}}
   \cdots z_{i_{(j-1)}}^{\alpha_{(j-1)}} z_{i_j}^{\alpha_{j}-1}
   z_{i_{(j+1)}}^{\alpha_{(j+1)}} \cdots z_{i_n}^{\alpha_{n}}},
\]
where, to ease the notation, we write
${z_{i_1}^{\alpha_{1}}z_{i_2}^{\alpha_{2}} \cdots
z_{i_n}^{\alpha_{n}}}$ instead of
$e_{{z_{i_1}^{\alpha_{1}}z_{i_2}^{\alpha_{2}} \cdots
z_{i_n}^{\alpha_{n}}}}$, and extend by linearity to any finite
linear combination of such elements, and prove that these
operators are continuous.

\begin{Pn}
$D_m$ is a well defined continuous linear operator $\widetilde
\bigS_{-1} \to \widetilde \bigS_{-1}$ and it holds that
\[
D_m(f\otimes g)=D_m(f)\otimes g+f\otimes D_m(g)
\]
for any $f,g \in \widetilde \bigS_{-1}$.
\end{Pn}
\begin{proof}[\indent Proof]
Let $f =\sum_{\alpha \in \widetilde \ell}f_\alpha e_\alpha \in
\widetilde\bigS_{-1}$. Then there exists $p \in \mathbb N$ such
that
\[
\sum_{\alpha \in \widetilde \ell}|f_\alpha|^2 (2\mathbb
N)^{-\alpha p}<\infty.
\]
For any $0\leq j \leq n$, let $r_j$ be defined by
\[
r_j\quad:\quad\{\alpha \in \widetilde \ell : |\alpha|=n\} \to
\{\alpha \in \widetilde \ell : |\alpha|=n+1\}
\] defined by
\[
r_j(z_{i_1}z_{i_2}\cdots z_{i_n})=z_{i_1}z_{i_2}\cdots
z_{i_j}z_{m}z_{i_{j+1}} \cdots z_{i_n}.
\]
Since $m$ is fixed, we do not write the dependence of $r_j$ on
$m$. Furthermore, we now allow $i_k=i_{k+1}$. Let $\beta \in
\widetilde \ell$. Then for any $\alpha\in \widetilde{\ell}$ and
for any $0 \leq j \leq |\alpha|$ such that $r_j(\alpha)=\beta$, we
have $|\alpha|+1=|\beta|$ and
\[
(2 \mathbb N)^\alpha=\prod_{l=1}^{|\alpha|}(2i^{(\alpha)}_l)=
(2m)^{-1}\prod_{l=1}^{|\beta|}(2i^{(\beta)}_l)=(2m)^{-1}(2\mathbb
N)^\beta.
\]
Moreover,
\[
\begin{split}
|\{(\alpha, j):\alpha \in \widetilde \ell, 0 \leq j \leq |\alpha|,
r_j(\alpha)=\beta\}|&=\\
&\hspace{-3cm}= |\{1 \leq k \leq |\beta| : \beta=z_{i_1}\cdots
z_{i_{|\beta|}},i_k=m\}| \leq |\beta|.
\end{split}
\]
Thus, denoting $\widetilde \ell_m=\{\beta \in \widetilde \ell:
\beta=z_{i_1}\cdots z_{i_{|\beta|}}, i_k=m \text { for some }k\}$
\begin{align*}
\|D_m f\|_q^2 &=\sum_{\alpha \in \widetilde \ell}
\left|\sum_{j=0}^{|\alpha|}f_{r_j(\alpha)}
\right|^2 (2\mathbb N)^{-\alpha q}\\
\displaybreak[1]&\leq \sum_{\alpha \in \widetilde
\ell}(|\alpha|+1)^2
\sum_{j=0}^{|\alpha|}\left|f_{r_j(\alpha)} \right|^2 (2\mathbb N)^{-\alpha q}\\
\displaybreak[1]&=\sum_{\beta \in \widetilde \ell_m
}\sum_{\{(\alpha, j):\alpha \in \widetilde \ell, 0 \leq j \leq
|\alpha|, r_j(\alpha)=\beta\}}
(|\alpha|+1)^2\left|f_{r_j(\alpha)} \right|^2 (2\mathbb N)^{-\alpha q}\\
\displaybreak[1]&=\sum_{\beta \in \widetilde \ell_m
}\sum_{\{(\alpha, j):\alpha \in \widetilde \ell, 0 \leq j \leq
|\alpha|, r_j(\alpha)=\beta\}}
|\beta|^2\left|f_{\beta} \right|^2 (2m)^q (2\mathbb N)^{-\beta q}\\
\displaybreak[1]&\leq\sum_{\beta \in \widetilde \ell_m }
|\beta|^3\left|f_{\beta} \right|^2 (2m)^q (2\mathbb N)^{-\beta q}.
\end{align*}
By induction it can be easily checked that for any $n \in \mathbb
N$, $2^{3(n-1)}\geq n^3$. Thus, for any $q \geq p +3$ and for any
$\beta \in \widetilde \ell_m$,
\[
(2m)^{-(q-p)}(2\mathbb N)^{(q-p)\beta}=(2m)^{-(q-p)} (2i_1 \cdots
2m \cdots  2i_{|\beta|})^{q-p} \geq 2^{3(|\beta|-1)}\geq
|\beta|^3.
\]
Therefore,
\[
|\beta|^3 (2m)^q (2\mathbb N)^{-\beta q} \leq (2m)^p(2\mathbb
N)^{-\beta p},
\]
and we obtain
\[
\|D_m f\|_q^2 \leq (2m)^p\|f\|_p^2.
\]
Hence, $D_m|_{\Gamma(\bigH_{-q})}:\Gamma(\bigH_{-p}) \to \widetilde \bigS_{-1}$ is bounded and therefore continuous.
Since as a strong dual of a reflexive Fr\'echet space $\widetilde \bigS_{-1}$ the inductive limit of the Hilbert spaces $\Gamma(\bigH_{-p})$, and by the universal property of inductive limits, we obtain that $D_m$ is continuous.\\

It is now easy to check that for any $f,g \in \widetilde
\bigS_{-1}$ which are finite linear combinations of the basis
$(e_\alpha)$, $D_m(f\otimes g)=D_m(f)\otimes g+f\otimes D_m(g)$.
By continuity it holds for any $f,g \in \widetilde \bigS_{-1}$.
\end{proof}

We recall that for a unital (associative) ring $R$ a pair
$(C,A)\in{ R}^{p\times N}\times{ R}^{N\times N}$ is called {\sl
observable} if there exists some $p \geq 0$ such that
\[
\begin{pmatrix}C & CA & C A^2
&\cdots& CA^{q-1}\end{pmatrix}
\]
is left invertible. If furthermore, we may choose $q=N$, then we
the pair $(C,A)$ is called {\sl strongly observable}.\\

In the following theorem and its proof we omit the symbol
$\otimes$ for simplicity.

\begin{Tm}
Let $(C,A)\in{\widetilde \bigS_{-1}}^{p\times N}\times{\widetilde
\bigS_{-1}}^{N\times N}$. If the pair $(E[C],E[A])$ is
observable, then the pair $(C,A)$ is observable. \label{tm:obs}
\end{Tm}

\begin{proof}[\indent Proof]
Let $q \geq 0$ be such that $\begin{pmatrix}E[C] & E[C]E[A] &
\cdots&E[ C]E[A^{q-1}]\end{pmatrix}$ is left invertible. We show
that for any $f \in (\widetilde \bigS_{-1})^{qN}$ such that
$$\begin{pmatrix}C & CA &
&\cdots& CA^{q-1}\end{pmatrix} f =0$$ it holds that $f=0$.\\

First, we note that for such $f$, $\begin{pmatrix}E[C] & E[C]E[A]
& \cdots&E[ CA^{q-1}]\end{pmatrix}E[f]=0$. Hence, $f_0=E[f]=0$.\\
Now,
\[
\begin{split}
0 &=(E  D_m)(\begin{pmatrix}C & CA
&\cdots& CA^{q-1}\end{pmatrix}f)\\
&=(E D_m)\begin{pmatrix}C & CA
&\cdots& CA^{q-1}\end{pmatrix}E[f]\\
&\quad+\begin{pmatrix}E[C] & E[C]E[A]
&\cdots& E[C]E[A^{q-1}]\end{pmatrix}(E  D_m)f\\
&=\begin{pmatrix}E[C] & E[C]E[A] &\cdots&
E[C]E[A^{q-1}]\end{pmatrix}f_{z_m}.
\end{split}
\]
implies $f_{z_m}=0$.

Furthermore, by a simple induction since there exist some
$\{U_k\}_{k<n}$ such that
\[
D_m^n(\begin{pmatrix}C & CA &\cdots& CA^{q-1}\end{pmatrix}f)
=\sum_{k<n}U_k D_m^nf +\begin{pmatrix}C & CA &\cdots&
CA^{q-1}\end{pmatrix}D_m^nf
\]
we conclude
\[
\begin{split}
0 &=(E  D_m^n)(\begin{pmatrix}C & CA
&\cdots& CA^{q-1}\end{pmatrix}f)\\
&=\sum_{k<n}E[U_k]E[D_m^nf] +\begin{pmatrix}E[C] & E[C]E[A]
&\cdots& E[C]E[A^{q-1}]\end{pmatrix}(E  D_m^n)f\\
&=\begin{pmatrix}E[C] & E[C]E[A] &\cdots&
E[C]E[A^{q-1}]\end{pmatrix}f_{z_m^n}.
\end{split}
\]
Thus, $f_{z_m^n}=0$, for any $n$ and $m$.\\
The next step is to show that $f_{z_lz_m}=0$. Since,
\[
\begin{split}
0 &=(E  D_lD_m)(\begin{pmatrix}C & CA
&\cdots& CA^{q-1}\end{pmatrix}f)\\
&=(E  D_l)(\begin{pmatrix}C & CA &\cdots&
CA^{q-1}\end{pmatrix}E[D_mf]
\\&\quad+(\begin{pmatrix}E[C] & E[C]E[A]
&\cdots& E[C]E[A^{q-1}]\end{pmatrix}E[D_lD_mf]
\\&\quad+(E  D_lD_m)(\begin{pmatrix}C & CA
&\cdots& CA^{q-1}\end{pmatrix}E[f]
\\&\quad+(E  D_m)(\begin{pmatrix}C & CA
&\cdots& CA^{q-1}\end{pmatrix}E[D_lf]\\
&=\begin{pmatrix}E[C] & E[C]E[A] &\cdots&
E[C]E[A^{q-1}]\end{pmatrix}f_{z_lz_m}
\end{split}
\]
we conclude that $f_{z_lz_m}=0$.\\

In the same manner it is easy to complete the proof and showing
that $f_\alpha=0$ for any $\alpha \in \widetilde \ell$.
\end{proof}

In the approach outlined here to non-commutative linear systems
we replaced the complex numbers by a non-commutative algebra with
a special topological structure. Other approaches are possible. We
mention in particular the work of Fliess \cite{MR51:583}. We also
mention \cite{MR2187742,MR2019348,MR1909375,MR2489365}.
Furthermore, using the setting developed in the present paper,
one can study non-commutative versions of stationary increments
stochastic processes and associated stochastic integrals in a way
similar to \cite{aal2,aal3}. This will presented in a future
publication. For related work on free stationary increments
stochastic processes, we refer to \cite{MR2540072,MR2770019}.\\

{\bf Acknowledgments:} It is a pleasure to thank Professor Marek
Bozejko for suggesting us to study the non-commutative version of
the Kondratiev space $\bigS_{-1}$ of stochastic distributions, and
for pointing our attention to the papers
\cite{MR2540072,MR2770019}. It is also a pleasure to thank
Professor Yuri Kondratiev for comments on this paper. Finally we
thank the referee for her/his comments.


\bibliographystyle{plain}

\def\cprime{$'$} \def\lfhook#1{\setbox0=\hbox{#1}{\ooalign{\hidewidth
  \lower1.5ex\hbox{'}\hidewidth\crcr\unhbox0}}} \def\cprime{$'$}
  \def\cprime{$'$} \def\cprime{$'$} \def\cprime{$'$} \def\cprime{$'$}

\end{document}